\documentclass[12pt]{amsart}
\usepackage{a4,amsfonts,amssymb,xy}

\xyoption{all}

\numberwithin{equation}{section}

\newtheorem{theo}{Theorem}[section]
\newtheorem{prop}[theo]{Proposition}
\newtheorem{lem}[theo]{Lemma}

\theoremstyle{definition}
\newtheorem{defi}[theo]{Definition}
\newtheorem{example}[theo]{Example}

\newtheorem{rem}[theo]{Remark}

\newcommand\N{{\mathbb{N}}}
\newcommand\R{{\mathbb{R}}}
\newcommand\C{{\mathbb{C}}}
\newcommand\uU{{\mathcal{U}}}
\newcommand\eE{{\mathcal{E}}}
\newcommand\kK{{\mathcal{K}}}
\newcommand\lL{{\mathcal{L}}}
\newcommand\tT{{\mathcal{T}}}

\newcommand\pf{\begin{proof}}
\newcommand\pfend{\end{proof}}

\begin{document}
\title{The coarse Baum--Connes conjecture and groupoids II}
\author{J. L. Tu}
\begin{abstract}
Given a (not necessarily discrete) proper metric space $M$ with bounded geometry, we define a groupoid $G(M)$. We show that the coarse Baum--Connes conjecture with coefficients, which states that the assembly map with coefficients for $G(M)$ is an isomorphism, is hereditary by taking closed subspaces.
\end{abstract}
\maketitle

\section*{Introduction}
Let $(X,d)$ be a metric space that we will suppose in this introduction to be uniformly locally finite for simplicity, i.e. $\forall R>0$, $\exists N\in\N$, $\forall x\in X$, $\# B(x,R)\le N$.

A subset $E$ of $X\times X$ is \emph{controlled} if $d_{|E}$ is bounded. Let $H$ be a separable, infinite dimensional Hilbert space. Let $C^*(X)$ be the closure of the algebra of operators $T\in \lL(\ell^2(X,H))$ whose support is controlled, such that every matrix element $T_{xy}\in \lL(H)$ is a compact operator.

For every real number $d>0$, let $P_d(X)$ be the space of probability measures on $X$ whose support have diameter $\le d$. Then the coarse Baum--Connes conjecture \cite{yu96} states that a certain assembly map $\lim_d K_*(P_d(X)) \to K(C^*(X))$ is an isomorphism.

This conjecture is known to be true in many cases \cite{yu00}, but not in general \cite{hls}.

In \cite{sty}, it was shown that $G(X)=\cup_{E\mbox{ controlled}} \bar{E}\subset \beta{X\times X}$ can be endowed with the structure of an \'etale, locally compact, $\sigma$-compact groupoid, and that the coarse Baum--Connes conjecture for $X$ is equivalent to the Baum--Connes conjecture for $G(X)$ with coefficients in $\ell^\infty(X,\kK)$.

In this paper, we extend the main result of \cite{sty} in two directions. First, we extend the construction to a large class of locally compact, proper metric spaces (that are not necessarily discrete). Secondly, we define a coarse Baum--Connes with coefficients: a natural way to do so is to require the groupoid $G(X)$ to satisfy the Baum--Connes conjecture with coefficients. We show that it is stable under taking closed subspaces. To that end, we prove that under quite general conditions on the locally compact groupoids $H\subset G$, the Baum--Connes conjecture with coefficients for $G$ implies the Baum--Connes conjecture with coefficients for $H$ (Theorems~\ref{thm:BCcoef1} and~\ref{thm:BCcoef}): this extends one of the main results in \cite{ceo}.

\section{General notations and conventions}
In a metric space, $B(a,R)$ (resp. $\tilde{B}(a,R)$) denotes the open ball (resp.
the closed ball) of center $a$ and radius $R$. More generally, if $A$ is a subspace then
$B(A,R)=\{x\vert\; d(x,A)<R\}$ and $\tilde{B}(A,R)=\{x\vert\; d(x,A)\le R\}$.

A metric space is said to be proper if all closed balls are compact.

If $G$ is a groupoid, we will denote by $G^{(0)}$ the space of units, and by $s$ and $r$ the source and the range maps.
For all $x,y\in G^{(0)}$, $G_x$, $G^y$ and $G_x^y$ denote $s^{-1}(x)$, $r^{-1}(y)$ and $G_x\cap G^y$. More generally, if $A,B\subset G^{(0)}$ then $G_A=s^{-1}(A)$, $G^B=r^{-1}(B)$ and $G_A^B=G_A\cap G^B$.

In particular, given a set $M$, $M\times M$ is endowed with the groupoid product $(x,y)(y,z)=(x,z)$ and inverse $(x,y)^{-1}=(y,x)$.

For all sets $A,B\subset M\times M$,
$A\circ B =\{(x,y)\in M\times M\vert\; \exists z\in M,\; (x,z)\in A\mbox{ and } (z,y)\in B\}$,
$A^{-1}=\{(y,x)\vert\; (x,y)\in A\}$,
$A_x=A\cap (M\times \{x\})$ and $A^x = A\cap (\{x\}\times M)$.
More generally, if $X\subset M$ then
$A_X=A\cap (M\times X)$ and $A^X = A\cap (X\times M)$. We will sometimes write $A\circ X$ instead of $A_X$.

Let $G$ be a groupoid. A right action of $G$ on a space $Z$ is given by a map $\sigma : Z\to G^{(0)}$ (the anchor map of the action) and a ``product'' $Z\times_{\sigma,r} G\to Z$, denoted by $(z,g)\mapsto zg$, satisfying the relations $z\sigma(z)=z$ and $(zg)h=z(gh)$ for all $(z,g,h)\in Z\times_{\sigma,r}G\times_{s,r}G$. A space endowed with an action of $G$ is called a $G$-space.

A continuous action is said to be proper if the map $Z\times_{\sigma,r}G\to Z\times Z$ defined by $(z,g)\mapsto (z,zg)$ is proper.
 
A space $Z$ endowed with an action of a groupoid $G$ is said to be $G$-compact (or cocompact) if $M/G$ is compact.

If a locally compact groupoid with Haar system acts properly on a locally compact space $Z$, then \cite{tu99} there exists a ``cutoff'' function $c:Z\to \R_+$ satisfying

\begin{itemize}
\item[(i)] $\forall x\in Z$, $\int_{g\in G^{\sigma(z)}} c(zg)\,\lambda^x(dg)=1$;
\item[(ii)] for every compact set $K\subset Z$, the set $\{(z,g)\in K\times G\vert\; c(zg)\ne 0\}$ is relatively compact.
\end{itemize}

\section{Uniform coarse structures and groupoids}
In this section, we associate to any LBG (see Proposition~\ref{prop:G'}) proper metric space $M$ a locally compact groupoid $G(M)$ (Definition~\ref{def:G(M)}). Most of the constructions below can be extended to spaces that are endowed with a uniform structure and a coarse structure which are compatible. However, we will deal most of the time with metric spaces, since spaces that one usually encounters are metrizable (see for instance Propositino~\ref{prop:exists-distance}).

We recall the following definition from general topology.

\begin{defi}
Let $M$ be a set. A uniform structure on $M$ is a nonempty
collection $\uU$ of subsets of $M\times M$ satisfying the following conditions:

\begin{itemize}
\item[(i)] For all $U\in \uU$,the diagonal $\Delta$ is a subset of $U$;
\item[(ii)] For all $U\in \uU$ and all $V\supset U$, we have $V\in \uU$;
\item[(iii)] For all $U,V\in \uU$, $U^{-1}\in \uU$ and $U\cap V\in\uU$;
\item[(iv)] For all $U\in \uU$, there exists $V\in \uU$ such that $V\circ V\subset U$.
\end{itemize}
\end{defi}

For instance, if $M$ is a metric space then $\uU$ consists of the subsets which contain
$\Delta_r=\{(x,y)\in M\times M\vert\; d(x,y)\le r\}$ for some $r$.

Given a uniform structure, there is a topology such that a subset $\Omega$ of $M$ is open if and only
for all $x\in \Omega$ there exists $U\in \uU$ satisfying the condition $U_x\subset \Omega$.
If a topological space $M$ is given, we call ``uniform structure on $M$''
a uniform structure which induces the topology on $M$.

A map
$f:M\to N$ between two uniform spaces is said to be uniformly continous if $(f\times f)^{-1}(V)\in \uU_M$
for all $V\in \uU_N$.

\begin{lem}
Let $\uU$ be a uniform structure on a topological space $M$. Given any neighborhood $W$ of the diagonal and $x\in M$,
there exists $V\in \uU$ and a neighborhood $\Omega$ of $x$ such that $V_\Omega\subset W$.
\end{lem}

\pf
Let $A$ be an open neighborhood of $x$ such that $A\times A\subset W$. Let $U\in \uU$ such that $U^x\subset A\times\{x\}$.
Let $V\in\uU$ such that $V^{-1}\circ V\subset U$. Since $V$ is a neighborhood of the diagonal, there exists an open
neighborhood $\Omega$ of $x$ such that $\Omega\times \{x\}\subset V$.

Let $(y,z)\in V_\Omega$, and let us prove that $(y,z)\in W$.
Since $(y,x)\in \Omega\times\{x\}\subset V\subset U$, we have
$y\in A$.

Since $(z,x)=(z,y)(y,x)\in V^{-1}\circ V\subset U$, we have $z\in A$. Therefore, $(y,z)\in A\times A\subset W$.
\pfend

If a group $\Gamma$ acts on a uniform space $M$, we will say that the uniform structure is $\Gamma$-invariant
if every $U\in \uU$ contains an element of $\uU$ which is $\Gamma$-invariant. For instance, a
$\Gamma$-invariant distance provides such a uniform structure.

\begin{prop}
Let $\Gamma$ be a locally compact group. Let $Y$ be a locally compact, $\Gamma$-compact proper $\Gamma$-space. Then
there is one and only one $\Gamma$-invariant uniform structure : a set $U$ belongs to $\uU$ if and
only if it contains a $\Gamma$-invariant neighborhood of the diagonal. As a consequence, if
$Z$ is any topological space with a $\Gamma$-invariant
uniform structure, then every continuous, $\Gamma$-invariant map $f:Y\to Z$ is uniformly continous.
\end{prop}

\pf
Let $W$ be a $\Gamma$-invariant neighborhood of the diagonal. We have to show that
$W\in \uU$.
Let $K\subset Y$ be a compact subset such that $K\Gamma = Y$.
By the preceding lemma, for all $x\in K$ there exists a neighborhood $\Omega_x$ of $x$
and $V_x\in \uU$ such that $V_x\cap(\Omega_x\times M)\subset W$. Let $x_1,\ldots,x_n\in K$
such that $K\subset\cup_i\Omega_{x_i}$. Let $U=\cap_i V_{x_i}$. Then $U\in \uU$, so there exists $U'\in \uU$
$\Gamma$-invariant contained in $U$. Since $U'\cap (K\times M)\subset W$, by invariance of $U'$ and of
$W$ we get $U'\subset W$.
\par\medskip
For the last statement, observe that if $U\in \uU_Z$ is $\Gamma$-invariant, then $(f\times f)^{-1}(U)$
is a $\Gamma$-invariant neighborhood of $\Delta_Y$, hence belongs to $\uU_Y$.
\pfend

\begin{defi} (Roe)
Let $M$ be a locally compact topological space. A coarse structure on $M$ is a collection $\eE$ of subsets of $M\times M$,
called \emph{entourages}, that have the following properties:
\begin{itemize}
\item[(a)] for any entourages $A$ and $B$, $A^{-1}$ and $A\circ B$ are entourages;
\item[(b)] any subset of an entourage is an entourage;
\item[(c)] every compact subset of $M\times M$ is an entourage.
\end{itemize}
\end{defi}

\begin{defi}
A uniform-coarse structure on a locally compact space $M$ is a pair $(\eE,\uU)$ consisting of
a coarse structure $\eE$, a uniform structure $\uU$, such that given $U\in \uU$ there exists
$V\subset U$ such that $V\in \uU\cap \eE$.
\end{defi}

For instance if $d$ is a proper distance of $M$ (meaning that every closed ball is compact) then,
with the coarse structure given by $E\in \eE \iff d_{|E}$ is bounded, and with the canonical
uniform structure, $M$ becomes a uniform-coarse space, which is \emph{proper} (in the sense that for
all $E\in \eE$, the projection maps $\bar{E}\to M$ are proper).

For most of the rest of the paper, we will deal with uniform-coarse structures which come from a metric.
Indeed, most locally compact spaces we will work with are metrizable (recall that a locally compact
space $X$ is metrizable if and only if $C_0(X)$ is separable, if and only if $X$ is second-countable, meaning
that its topology has a countable basis). Moreover, we have the following proposition.

\begin{prop}\label{prop:exists-distance}
Let $\Gamma$ be a locally compact group acting properly on a locally compact space $Y$ such that
$Y/\Gamma$ is compact. There is one and only one uniform-coarse structure on $Y$ which is proper and
$\Gamma$-invariant (a coarse structure is $\Gamma$-invariant if every entourage is
contained in a $\Gamma$-invariant one): entourages consist of sets $E\subset Y\times Y$
which are $\Gamma$-relatively compact (i.e. contained in a $\Gamma$
-invariant, $\Gamma$-compact set).
Moreover, if $\Gamma$ is discrete then there exists
a $\Gamma$-invariant (proper) distance on $Y$ which induces the above-mentioned uniform-coarse
structure.
\end{prop}

\pf
To show the existence part in the first assertion, we have to prove that every $\Gamma$-invariant neighborhood of $\Delta$
contains a $\Gamma$-compact $\Gamma$-invariant neighborhood of $\Delta$. This follows from the fact that
$\Gamma$-invariant open sets in $Y\times Y$ correspond to open subsets of $(Y\times Y)/\Gamma$,
and that $(Y\times Y)/\Gamma$ is locally compact.

To show uniqueness, let $\eE'=\{E\subset Y\times Y\;
\Gamma-\mbox{relatively compact}\}$.

Let $(\uU,\eE)$ be a uniform-coarse proper $\Gamma$-invariant
structure. Since $Y$ is $\Gamma$-compact, we have $\eE\subset \eE'$.

Conversely, since $\eE$ contains all compact subsets (by definition
of a coarse structure) and is $\Gamma$-invariant, we have $\eE'\subset \eE$.

Let us show the last assertion. Let $d$ be a distance on $Y$. After replacing $d(x,y)$
by $d(x,y)+|\varphi(x)-\varphi(y)|$ where $\varphi:Y\to \R$ is a proper continuous function,
we may assume that $d$ is a proper distance. Choose $y_0\in Y$ and $R>0$ such that $K\Gamma=Y$,
where $K$ is the closed ball $\tilde{B}(y_0,R)$. For all $n\ge 1$, let
$(c_{n,i})_{1\le i\le i_n}$ be a finite family of functions $c_{n,i}\in C_c(Y)_+$ such that
$\mbox{diam}\,(\mbox{supp}\, c_{n,i}) \le 2^{-n}$ and $K\subset \cup_{i=1}^{i_n} c_{n,i}^{-1}(\R^*)$
and $\sup_{y\in Y} \sum_\gamma c_{n,i}(y\gamma) \le 2^{-n-i}$. Consider
$$d_1(y,y')=\sum_{n,i}\int |c_{n,i}(y\gamma)-c_{n,i}(y'\gamma)|\,d\gamma.$$
Then $d_1$ is a $\Gamma$-invariant distance. To see this, the only non-obvious part is to check that
if $d_1(y,y')=0$ then $y=y'$. Let $L=\tilde{B}(y_0,R+1)$. Let $F$ be the closure of
$\{\gamma\in\Gamma\vert\; L\gamma\cap L\ne\emptyset\}$. Then $F$ is finite. For all $n$, there
exists $\gamma_n\in\Gamma$ such that $c_{n,i}(y\gamma_n)\ne 0$. Since $c_{n,i}(y'\gamma_n)\ne 0$,
it follows that $d(y\gamma_n,y'\gamma_n)\le 2^{-n}$. Since $F$ is finite, there exists $\gamma\in F$
such that $d(y\gamma,y'\gamma)=0$, so that $y=y'$.
\par\medskip

Let $c\in C_c(Y)_+$ such that $\sum_\gamma c(y\gamma)=1$ for all $y$. Let $P_r(\Gamma)$ the simplicial
set such that simplices consist of subsets of $\Gamma$ of diameter $\le r$. Then
$\mu:y\mapsto \sum_\gamma c(y\gamma)\delta_\gamma$ determines a $\Gamma$-equivariant map
from $Y\to P_r(\Gamma)$ for some $r$, thus determines a function $d_2:Y\times Y\to \R_+$ which
satisfies all the properties of a proper $\Gamma$-invariant distance except perhaps for
the separation axiom. Then $d_1+d_2$ is a $\Gamma$-invariant distance on $Y$.
\pfend

\begin{defi}
Let $X$ be a metric space. We say that $X$ is ULF (uniformly locally finite) if for all $R>0$,
$\sup_{x\in X}\# B(x,R) < +\infty$.
\end{defi}

\begin{defi}
A metric space $X$ is $\delta$-separated (resp. strictly $\delta$-separated) if $d(x,y)\ge \delta$ (resp. $d(x,y)>\delta$) for all $x\ne y \in X$.
\end{defi}

\begin{defi}
Let $M$ be a metric space. A subset $X$ is said to be $\varepsilon$-dense (resp. strictly $\varepsilon$-dense) if for all $m\in M$,
$d(m,X)<\varepsilon$.
\end{defi}

\begin{defi}
A metric space $M$ is said to have bounded geometry if for all $\varepsilon>0$ there exists a subspace
$X$ which is $\varepsilon$-dense and ULF.
\end{defi}

\begin{example}
Let $M$ be a compact Riemannian manifold. Then the universal cover of $M$ has bounded geometry.
\end{example}

\pf
Let $\Gamma$ be the fundamental group of $M$. Let $\pi:\tilde{M}\to M$ be the natural projection.
Let $X\subset M$ finite such that $\cup_{x\in X}B(x,\varepsilon)=M$.
Let $\tilde{X}=\pi^{-1}(X)$. Then $\tilde{X}$ is $\varepsilon$-dense. Moreover, it is a finite union of
$\Gamma$-orbits, thus it is ULF.
\pfend

\begin{lem}
Let $M$ be a bounded geometry metric space. Then for all $R>0$ and all $\varepsilon>0$, there exists
$n\in \N$ such that for every nonempty
subset $A$ of $M$ of diameter $\le R$, there exist $a_1,\ldots,a_n\in A$
such that $A\subset\cup_{i=1}^n\tilde{B}(a_i,\varepsilon)$.
\end{lem}

\pf
Let $X\subset M$ ULF and $\varepsilon/2$-dense. Let $Y=\{x\in X\vert\;
d(x,A)< \varepsilon/2\}$. Let $n$ such that for all $Z\subset X$ of diameter $\le R+\varepsilon$,
we have $\# Z\le n$.

For all $x\in Y$, choose $f(x)\in A$ such that $d(x,f(x))\le \varepsilon/2$. Let $B=f(Y)$.
Then $\#B\le n$, and $A\subset\cup_{b\in B}\tilde{B}(b,\varepsilon)$.
\pfend

\begin{lem}
Let $N$ be an integer. Let $\Delta$ a graph such that each vertex has at most $N-1$ neighbors. Then one can color the vertices using at most $N$ colors, so that two neighboring vertices have different colors.
\end{lem}

\pf
We may assume that the graph is connected, hence countable. Label the vertices as $\{x_0,x_1,\ldots\}$. Suppose colors have been attributed to $x_0,\ldots,x_n$. Let $A_n$ be the set of colors of those $x_i$'s ($i\le n$) which are adjacent to $x_{n+1}$. Since $\# A_n\le N-1$, one can give to $x_{n+1}$ a color which does not belong to $A_n$.
\pfend

\begin{prop}\label{prop:N-R-separated}
Let $R>0$ and $N\in \N$. Let $X$ be a metric space such that every ball of radius $R$ has at most $N$ elements. Then there exists a decomposition $X=X_1\cup\cdots\cup X_N$ into $N$ strictly $R$-separated spaces.
\end{prop}

\pf
Apply the preceding lemma to the graph whose vertex set is $X$, such that $(x,y)$ is an edge if and only if $x\ne y$ and $d(x,y)\le R$.
\pfend

We denote by $UC_b(M)$ the algebra of bounded, uniformly continuous functions on $M$. This is a (usually non-separable) abelian $C^*$-algebra. Let $\beta_uM$ be its spectrum. Note that $M$ is an open dense subset of the compact set $\beta_uM$.

The following property will be needed later:

\begin{lem}\label{lem:UCbMF}
Let $F$ be a closed subset of a locally compact metric space $M$. Then the restriction map $UC_b(M)\to UC_b(F)$ is surjective (and thus $\beta_uF$ can be identified with the closure of $F$ in $\beta_uM$).
\end{lem}

\pf
Let $f\in UC_b(F)$. Let $r(x)=d(x,F)$. Define
$$g(x)=\left\{\begin{array}{ll}
f(x) & \mbox{if } x\in F\\
\frac{1}{r(x)} \int_{r(x)}^{2r(x)} \inf_{\tilde{B}(x,t)\cap F}f\,dt & \mbox{otherwise.}
\end{array}
\right.
$$
We show that $g$ is uniformly continuous. After translating and rescaling, we may assume that $0\le f\le 1$. Let $\varepsilon\in (0,1)$. There exists $\eta\in (0,\varepsilon)$ such that $d(x,y)<\eta\implies |f(x)-f(y)|<\varepsilon$.
Let $x,y\in M$ such that $d(x,y)\le\eta^2/100$.

1st case : suppose $r(x)\ge\eta/5$. Let $h(z)=\int_{r(z)}^{2r(z)} \inf_{\tilde{B}(z,t)\cap F}f\,dt$.
\begin{eqnarray*}
h(x) &=& \int{r(x)}^{r(x)+2d(x,y)}  \inf_{\tilde{B}(x,t)\cap F}f\,dt
+ \int_{r(x)+2d(x,y)}^{2r(x)}  \inf_{\tilde{B}(x,t)\cap F}f\,dt\\
&\le& 2d(x,y)+ \int_{r(x)+d(x,y)}^{2r(x)-d(x,y)}  \inf_{\tilde{B}(x,t+d(x,y))\cap F}f\,dt\\
&\le& 2d(x,y)+ \int_{r(x)+d(x,y)}^{2r(x)-d(x,y)}  \inf_{\tilde{B}(y,t)\cap F}f\,dt \\
&\le& 3d(x,y)+h(y)
\end{eqnarray*}
and similarly $h(y)\le 3d(x,y)+h(x)$, so $|h(x)-h(y)|\le 3 d(x,y)$.

\begin{eqnarray*}
|g(x)-g(y)| &=& \left| \frac{h(x)r(y)-h(y)r(x)}{r(x)r(y)} \right | \\
&\le& \frac{|h(x)-h(y)| r(x) + |h(x)| \,|r(x)-r(y)|}{r(x)r(y)}\\
&\le& \frac{3d(x,y)}{r(x)-d(x,y)}+\frac{d(x,y)}{r(x)-d(x,y)}\\
&=& \frac{4d(x,y)}{r(x)-d(x,y)} \le \frac{4\eta^2/100}{\eta/5-\eta^2/100} \le\varepsilon.
\end{eqnarray*}

2nd case: if $r(y)\ge \eta/5$ then similarly $|g(x)-g(y)|\le\varepsilon$.

3rd case: suppose that $r(x)$ ,$r(y)< \eta/5$. We treat the case $r(x)>0$ and $r(y)>0$, the case $r(x)=0$ or $r(y)=0$ being similar.

$\forall t\in [r(x),2r(x)]$, $\forall s\in [r(y),2r(y)]$, $\forall u\in\tilde{B}(x,t)\cap F$, $\forall v\in \tilde{B}(y,s)\cap F$,
$d(u,v)\le d(u,x)+d(x,y)+d(y,v)\le t+\eta^2/100+s\le 4\eta/5+\eta^2/100\le\eta$, so $|g(x)-g(y)|\le\varepsilon$.

This completes the proof that $g$ is uniformly continuous. It is obviously bounded and extends $f$, so $f$ is in the image of $UC_b(M)\to UC_b(F)$.
\pfend

\begin{prop}\label{prop:bounded-geometry-ULF}
Let $M$ be a locally compact metric space, and let $\delta>0$.
The following are equivalent:
\begin{itemize}
\item[(i)] $M$ has bounded geometry;
\item[(ii)] $\beta_uM=\cup_{X\;\mathrm{ULF}, X\subset M}\bar{X}$;
\item[(iii)] $\beta_uM=\cup_{X\;\mathrm{ULF},\;\delta-\mathrm{sep.},\;X\subset M}\bar{X}$.
\end{itemize}
\end{prop}

\pf
(iii)$\implies$(ii) is obvious. To show the converse, we use the fact that if $X$ is ULF then it is a finite union of $\delta$-separated spaces (see Proposition~\ref{prop:N-R-separated}
).
\par\medskip

Let us show $(i)\implies (ii)$. Let $\alpha\in \beta_uM$.
Choose $X_0\subset M$ ULF, 1-dense. From the preceding lemma,
there exists $n_1\in\N$, and
for all $x\in X_0$ there exist $a_1^{(1)}(x),\ldots,a_{n_1}^{(1)}(x)\in\tilde{B}(x,1)$
such that $\tilde{B}(x,1)\subset\cup_i\tilde{B}(a_i^{(1)}(x),1/2)$. Let $B_{1,i}(x)=
\tilde{B}(a_i^{(1)}(x),1/2) \cap\tilde{B}(x,1)$ and
$A_i=\cup_{x\in X_0} B_{1,i}(x)$,
then $M=\cup_{i=1}^{n_1}A_i$, so there exists
$i_1$ such that $\alpha\in \bar{A}_{i_1}$. Let $B_1(x)=B_{1,i_1}(x)$.
Continuing in the same way, we cover $B_1(x)$ by balls of radius $1/4$, etc. and thus we get
$B_k(x)\subset B_{k-1}(x)\subset\cdots\subset \tilde{B}(x,1)$ compact such that
$\alpha\in \overline{\cup_{x\in X_0}B_k(x)}$ and $B_k(x)$ is of diameter $\le 2^{1-k}$.
Let $Y=\cup_{x\in X_0}\cap_{k\ge 1} B_k(x)$. Obviously, $Y$ is ULF. We want to show that $\alpha
\in \bar{Y}$. If this was not the case, there would exist $f\in UC_b(M)$ such that $f(\alpha)=1$
and $f_{|Y}=0$. By uniform continuity, we get $f\le 1/2$ on $\cup_{x\in X_0} B_k(x)$ for
$k$ large enough, and by continuity of $f$ at $\alpha$ we get $f(\alpha)\le 1/2$. Contradiction.
\par\medskip

Let us show $(ii)\implies (i)$. Suppose that for some $R>0$, ULF subsets are not $R$-dense.
For all $X\subset M$ ULF, denote by $f_X$ the function
$f_X(x)=\left(1-\frac{d(x,X)}{R}\right)_+$. Then
$f\in UC_b(M)$. We identify $f$ with a continuous function on $\beta_uM$. Let $F_X=\{x\in \beta_uM\vert\; f_X(x)= 0\}$. Since $X$ is not $R$-dense,
$F_X$ is a nonempty closed subset of $\beta_uM$. Moreover, if $X\subset Y$ then $F_X\supset F_Y$.
By compactness of $\beta_uM$, there exists $\alpha\in\beta_uM$ such that $\alpha\in F_X$ for all
$X\subset M$ ULF. Since $f_X(\alpha)=0$ and $f_X=1$ on $X$, we have $\alpha\notin \bar{X}$ for all
$X\subset M$ ULF.
\pfend

From now on, $(M,d)$ denotes  a bounded geometry locally compact proper metric space.
To understand better the topology of $\beta_uM$, we describe a basis of neighborhoods for each point
of $\beta_uM$.

\begin{prop}\label{prop:basis-neighborhoods}
Let $(M,d)$ be a bounded geometry, locally compact proper metric space. Let $\alpha\in \beta_uM$. Choose $X\subset M$ $\delta$-separated such that $\alpha\in \bar{X}$.
For each $Y\subset X$ such that $\alpha\in \bar{Y}$ and each $\varepsilon>0$, let $N_{Y,\varepsilon}=
\overline{B(Y,\varepsilon)}$. Then the $N_{Y,\varepsilon}$ constitute a basis of neighborhoods
of $\alpha$.
\end{prop}

\pf
Let $W$ be a neighborhood of $\alpha$. There exists $f\in UC_b(M)$ such that
$f(\alpha)=1$ and $f$ is supported in $W$. Let $\varepsilon>0$ such that
$d(x,y)\le \varepsilon$ implies $|f(x)-f(y)|\le 1/3$. Let $Y=\{x\in X\vert\; f(x)\ge 2/3\}$.
Then $\alpha\in \bar{Y}$, and for all $x\in B(Y,\varepsilon)$ we have $f(x)\ge 1/3$, so $f\ge 1/3$
on $N_{Y,\varepsilon}$, which implies that $N_{Y,\varepsilon}\subset W$.

Conversely, if $Y$ and $\varepsilon$ are as in the proposition, let
$f(x)=(1-d(x,Y))_+$. Since $f\in UC_b(M)$, $f$ extends to a continuous function
$h$ on $\beta_uM$. Since $h(\alpha)=1$, $U=h^{-1}((1-\varepsilon,1])$ is an open
neighborhood of $\alpha$.
Moreover, $U\cap M=B(Y,\varepsilon)$, so for all $\beta\in U$ and for every open neighborhood $V$ of
$\beta$, we have $V\cap B(Y,\varepsilon)=(V\cap U)\cap M\ne\emptyset$, which shows that
$\beta\in \overline{B(Y,\varepsilon)}=N_{Y,\varepsilon}$, for all $\beta\in U$, i.e. $U\subset N_{Y,\varepsilon}$.
This shows that $N_{Y,\varepsilon}$ is a neighborhood of $\alpha$.
\pfend

Our goal is now to define a groupoid associated to a locally compact proper bounded geometry metric space $(M,d)$.

We need some preliminaries.

Let $A_M$ be the abelian $C^*$-algebra consisting of $f\in UC_b(M\times M)$ such that
for all $\varepsilon>0$ there exists an entourage $E\in \eE$ such that
$|f|\le \varepsilon$ outside $E$. We define $G(M)$ as the spectrum of $A_M$. Since $C_0(M\times M)$
is an essential ideal of $A_M$, $M\times M$ is a dense open subset of $G(M)$.
Our goal is to show that the groupoid product $(x,y)(y,z)=(x,z)$ on $M\times M$ extends by continuity to
$G(M)$.

\begin{lem}
Let $X$ be a closed subspace of $M$. The restriction map $A_M\to A_X$ is surjective, and identifies
$G(X)$ with the closure of $X\times X$ in $G(M)$.
\end{lem}

\pf
Given an entourage $E$, let $A_{M,E}$ be the set of all $f\in A_M$ such that $f=0$ outside $E$.
We will write $A_{X,E}$ instead of $A_{X,E\cap(X\times X)}$ for simplicity.
Clearly, the union of all $A_{M,E}$ is dense in $A_M$, so it suffices to show that every $f\in A_{X,E}$ is the restriction of some element in $A_M$.
Indeed, since $f\in UC_b(X\times X)$, we already know (Lemma~\ref{lem:UCbMF}) that $f$ is the restriction of some
function $g\in UC_b(M\times M)$. Let $h(z)=g(z)(1-d(z,E))_+$. Then $h_{|X\times X}=f$ and $h\in A_M$.
This shows the first assertion, thus $G(X)$ is a subspace of $G(M)$. Since $X\times X$ is dense
in $G(X)$, $G(X)$ is the closure of $X\times X$ in $G(M)$.
\pfend

\begin{lem}\label{lem:subspace-bounded-geometry}
Let $M$ and $N$ be two bounded geometry metric spaces, and $A\subset M$. Then $M\times N$ and $A$
have bounded geometry.
\end{lem}

\pf
The first assertion is clear. Let us prove the second one. Let $R>0$. Choose an ULF subspace $X$ of
$M$ which is $R/3$-dense. Let $f:X\to A$ a map such that $d(f(x),x)\le 2 d(x,A)$ for all $x$.
We show that $Y=\{f(x)\vert\;x\in X,\;d(x,A)\le R/3\}$ is $R$-dense and ULF.

For all $a\in A$, there exists $x\in X$ such that $d(a,x)\le R/3$. Then $d(a,f(x))\le d(a,x)+d(x,f(x))
\le R/3+2d(x,A)\le R$, so $Y$ is $R$-dense.

$B_Y(f(x),S)\subset Y\cap B(x,S+2R/3)\subset
\{f(x')\vert\; d(x,x')\le 2R/3+S+2R/3\} =f(B(x,S+4R/3)))$, so $Y$ is ULF.
\pfend

\begin{lem}\label{lem:union-of-bar-E}
Let $M$ be a bounded geometry, locally compact proper metric space. $G(M)=\cup\bar{E}$, where $E$ runs over all entourages.
\end{lem}

\pf
Let $\alpha\in G(M)$. There exists $f\in A_M$ such that $f(\alpha)=1$. Let $E=\{z\vert\;f(z)>1/2\}$,
then $E$ is an entourage. Moreover, $\alpha\in \overline{M\times M}
=\bar{E}\cup\overline{\{z\vert\;f(z)\le 1/2\}}$. Since $\alpha\notin \overline{\{z\vert\;f(z)\le 1/2\}}$,
it follows that $\alpha\in \bar{E}$.
\pfend

\begin{prop}
Let $M$ be a bounded geometry, locally compact proper metric space. Then $G(M)=\cup_{X\subset M\;\mathrm{ULF}} G(X)$.
\end{prop}

\pf
Let $\alpha\in G(M)$. According to Lemma~\ref{lem:union-of-bar-E}, there exists a closed entourage
$E$ such that $\alpha\in \bar{E}$. We want to show that $\alpha$ is in
$\overline{X\times X}$ for some $X\subset M$ which is ULF. First, $E$ has bounded geometry since it is
a subspace of $M\times M$ (see Lemma~\ref{lem:subspace-bounded-geometry}).

According to Proposition~\ref{prop:bounded-geometry-ULF}, there exists $Y\subset E$ ULF such that
$\alpha\in\bar{Y}$. Let $X=\mbox{pr}_1(Y)\cup\mbox{pr}_2(Y)$. Since $\alpha\in \overline{X\times X}$, it
just remains to prove that $X$ is ULF. Let us show for instance that $X_1:=\mbox{pr}_1(X)$ is ULF.
Let $S=\sup_E d$. If $(x,y)$, $(x',y')\in Y$ satisfy $d(x,x')\le R$ then
$d(x,x')+d(y,y')\le R+R+2S=2R+2S$, so $B_{X_1}(x,R)\subset
\mbox{pr}_1(B_Y(y,2R+2S))$.
\pfend

\begin{lem}\label{lem:g1gn}
Let $X$ be a ULF metric space. Let $g_1,\ldots,g_n\in G(X)$. Then there exists $\delta>0$
and $X'\subset X$ $\delta$-separated such that $g_1,\ldots,g_n\in G(X')$.
\end{lem}

\pf
We use an induction over $n$. For $n=0$ there is nothing to prove.
Suppose that there exists $X'\subset X$ $\delta$-separated such that
$g_i\in G(X')$ for all $i<n$. Let $N$ such that balls of radius $\delta$ have at most $N$ elements.
Choose $\varepsilon\in (0,\delta/N)$. We define an equivalence relation $x\sim y$ on $X$ if there exists
$k$ and $x=x_0,\ldots,x_k=y$ such that $d(x_i,x_{i+1})\le \varepsilon$. Let $X_i$ ($i\in I$)
be the equivalence classes.
We have $\mbox{diam}\,(X_i)<\delta$, and $d(X_i,X_j)>\varepsilon$ if $i\ne j$.

Let $J=\{i\in I\vert\; X_i\cap X'\ne \emptyset\}$. As $X'$ is $\delta$-separated, for all $i\in J$
there exists $x_i$ such that $X_i\cap X'=\{x_i\}$.

Let $f_1(x,y)=\min(d(x,X'),1)$, $f_2(x,y)=\min(d(y,X'),1)$ and $f=\max(f_1,f_2)$.
Since $f_1$ and $f_2$ are uniformly continuous
and bounded, they are multipliers of $A_X$, thus they extend to continuous and bounded functions
$h_1$ and $h_2$ on $G(X)$. Let $h=\max(h_1,h_2)$.

1st case:
Suppose that $h(g_n)=0$. If $g_n\notin \overline{X'\times X'}$, then there exists $\varphi:X\times X\to\R$
uniformly continuous such that $\varphi(g_n)=1$ and $\varphi=0$ on $X'\times X'$. Since $\varphi$
is uniformly continuous, there exists $\eta\in (0,1)$ such that $f(x,y)<\eta$ implies
$\varphi(x,y)\le 1/2$. As a consequence, $g_n\notin\overline{\{(x,y)\vert\;f(x,y)<\eta\}}$,
so that $g_n\in\overline{\{(x,y)\vert\; f(x,y)\ge\eta\}}$. By continuity of $h$ we get $h(g_n)\ge \eta$.
Contradiction. This shown that $g_n\in\overline{X'\times X'}$.

2nd case:
$h(g_n)>0$. Suppose for definiteness that $h_1(g_n)>0$. Let $\eta\in (0,h_1(g_n))$.
Then $g_n\notin \overline{\{(x,y)\vert\; f_1(x,y)\le\eta\}}$, so $g_n\in
\overline{\{(x,y)\vert\; f_1(x,y)\ge\eta\}}$. For all $i\in I$, let
$x_{i,1},\ldots,x_{i,n_i}$ be the elements of $X_i$ such that $d(x_{i,\lambda},X')\ge\eta$.
We have $n_i\le N$ for all $i$. Let
$Y_\lambda =\{x_{i,\lambda}\vert\; i\in I\}$. Since $g_n\in \overline{(\cup_\lambda Y_\lambda)\times X}$,
there exists $\lambda$ such that
$g_n\in  Y_\lambda\times X$. After replacing $\delta$ by $\min(\delta,\eta,\varepsilon)$
and $X'$ by $X'\cup Y_\lambda$, we can assume that $g_n\in \overline{X'\times X}$,
thus that $h_1(g_n)=0$. Similarly, we can assume that $h_2(g_n)=0$, so we are reduced to the first case
treated above.
\pfend

Let us now define the product on the groupoid $G(M)$.
First, the source map $s(x,y)=y$ for the pair groupoid $M\times M$ defines a map
$UC_b(M)\to UC_b(M\times M)$, thus a map $\beta_u(M\times M)\to\beta_uM$.
In particular, $s$ extends continuously to a map $s:G(M)\to \beta_uM$.

If $(g,h)\in G(M)^2$ is a composable pair, from Lemma~\ref{lem:g1gn} there exists $X$
ULF $\delta$-separated such that $(g,h)\in G(X)$. Since $G(X)$ is a groupoid \cite{sty},
we can define the product in the groupoid $G(X)\subset G(M)$. Let us show that the product
does not depend on the choice of $X$. Suppose that $X'$ and $X''$ are $\delta$-separated
and that $g,h\in G(X')\cap G(X'')$.

\begin{lem}\label{lem:X'X''}
Let $\alpha\in \bar{X}'\cap\bar{X}''$.
For all $\varepsilon>0$, let
\begin{eqnarray*}
X'_\varepsilon &=& \{x\in X'\vert\; d(x,X'')\le\varepsilon\}\\
X''_\varepsilon &=& \{x\in X''\vert\; d(x,X')\le\varepsilon\}.
\end{eqnarray*}
Then $\alpha\in \bar{X}'_\varepsilon \cap\bar{X}''_\varepsilon $.
\end{lem}

\pf
We can assume $0<\varepsilon<\min(1,\delta/2)$.
Since $f(x)=\max(d(x,X''),1)$ is uniformly continous, it extends to $h\in C(\beta_uM)$.
Since $h=0$ on $X''$ and $\alpha\in\overline{X}''$, we have $h(\alpha)=0$, so
$\alpha\notin \overline{f^{-1}([\varepsilon,1])}$. Thus, $\alpha\in \overline{X'_\varepsilon}$
and similarly, $\alpha\in \overline{X''_\varepsilon}$.
\pfend

Applying Lemma~\ref{lem:X'X''} to $X\times X$ and $X'\times X'$, we see that
$g,h\in G(X'_\varepsilon)\cap G(X''_\varepsilon)$. Now, for $\varepsilon<\delta/2$,
there exists a unique bijection $\varphi_\varepsilon:X'_\varepsilon\to X''_\varepsilon$
such that $d(x,\varphi_\varepsilon(x))\le\varepsilon$ for all $x\in X'_\varepsilon$.
This induces an isomorphism of groupoids, again denoted by $\varphi_\varepsilon$.
Let $\gamma'$ (resp. $\gamma''$) be the product of $g$ and $h$ computed in
$G(X')$ (resp. $G(X'')$). Since $\gamma''=\varphi_\varepsilon(\gamma')$, it suffices to
show that $\varphi_\varepsilon(g)=g$ and $\varphi_\varepsilon(h)=h$. Let us show for instance
$\varphi_\varepsilon(g)=g$. Note that $\varphi_\varepsilon(g)$ does not depend on $\varepsilon$.
If $\varphi_\varepsilon(g)\ne g$ then there exists $h\in UC_b(M\times M)$ such that
$h(\varphi_\varepsilon(g))=0$ and $h(g)=1$. Let $\varepsilon\in (0,\delta/2)$ such that
$d(\gamma_1,\gamma_2)\le \varepsilon \implies |h(\gamma_1)-h(\gamma_2)|\le 1/2$. Then
$ |h(\gamma)-h(\varphi_\varepsilon(\gamma))|\le 1/2$ for all $\gamma\in G(X'_\varepsilon)$,
so $|h(g)=h(\varphi_\varepsilon(g))|\le 1/2$. Impossible.
This completes the proof that the product in $G(M)$ is well-defined.

Let us show that the product is continuous. Suppose that $g,h\in G(X)$ are composable, where $X$
is $\delta$-separated. We want to show that if $W$ is a neighborhood of $gh$ then there exists
a neighborhood $U$ of $(g,h)$ such that for all composable $(g',h')\in U$ we have
$g'h'\in W$. Let $\varphi\in UC_b(M\times M)$ such that $\varphi(gh)=1$ and $\varphi$ is
supported in $W$. There exists $\eta\in (0,\delta/2)$ such that $d(x,x')\le \eta$ and $d(y,y')\le\eta$
imply $|\varphi(x,y)-\varphi(x',y')|\le 1/3$.

Choose an entourage $E$ such that $g,h\in \bar{E}_X$, where
$E_X=E\cap(X\times X)$. From \cite{sty}, there exist $F_1,F_2\subset E_X$ such that
the source and range maps are injective on $F_1$ and $F_2$, $(g,h)\in F_1\times_X F_2$, and $\varphi(g'h')\ge 2/3$
for all $(g,'h')\in F_1\times_X F_2$. Let
$F'_i=B(F_i,\eta)$, then $\bar{F}'_i$ are neighborhoods of $g$ and $h$ respectively such that
$\varphi(g''h'')\ge 1/3$ for all $(g'',h'')\in F'_1\times_M F'_2$.
By continuity, $\varphi(g'',h'')\ge 1/3$ for $(g'',h'')$ in a neighborhood of $(g,h)$, which proves that
$g''h''\in W$.
\par\medskip
This proves that the product in $G(M)$ is continous. The fact that the inverse map $g\mapsto g^{-1}$
is even simpler.
\par\medskip

The groupoid $G(M)$ is $\sigma$-compact: indeed, $G(M)$ is the union of
$$\overline
{\{(x,y)\in M\times M\vert\; d(x,y)\le n\}}.$$
\par\medskip
There exists a Haar system on $G(M)$. To see this, we need a
\begin{lem}
There exists a measure $\mu$ on $M$ such that for all $R>0$,
\begin{itemize}
\item[(i)] $\sup_{x\in M}\mu(B(x,R)) < \infty$;
\item[(ii)] $\inf_{x\in M}\mu(B(x,R)) > 0$.
\end{itemize}
\end{lem}

\pf
For all $n\ge 1$, let $X_n\subset M$ ULF and $1/n$-dense. Let $a_n(R)=\sup_{x\in X_n}\#B_{X_n}(x,R)$,
$\mu_n=\sum_{x\in X_n}\delta_x$, $c_n=2^{-n}(1+a_n(n))^{-1}$ and $\mu=\sum_{n\ge 1} c_n\mu_n$.

Let us prove (i).
For all $x\in M$ and $n\ge 1$, there exists $y\in X_n$ such that $d(x,y)\le 1$.
Since $\mu_n(B(x,R))\le \mu_n(y,R+1)=\# B_{X_n}(y,R+1)\le a_n(R+1)$,
we have $\mu(B(x,R))\le \sum_{n=1}^\infty 2^{-n}a_n(R+1)(1+a_n(n))^{-1}<\infty$.

Let us prove (ii).
Let $n>1/R$. Then $\mu(B(x,R))\ge c_n\mu_n(B(x,R))\ge c_n>0$.
\pfend

\begin{rem}
In fact, the existence of a measure satisfying properties (i) and (ii) above is equivalent to
the fact that $M$ has bounded geometry.
\end{rem}

We now define the Haar system as follows.

The $C(\beta_uM)$-linear map $f\in C_c(G(M))\mapsto \varphi\in C(\beta_uM)=UC_b(M)$ defined by
$$\varphi(x)=\int_M f(x,y)\,d\mu(y)$$
defines a Haar system $(\lambda^x)_{x\in\beta_uM}$. Indeed, the fact that $\varphi$ is well-defined is a consequence of (i), and the fact that $\lambda^x$ has support $G(M)^x$ is a consequence of (ii).

Now, we generalize the definition of $G(M)$ to metric spaces that do not necessarily have bounded geometry.

\begin{defi}
Let $M$ be a metric space. We denote by $\eE'_M$ (or by $\eE'$ if there is no ambiguity) the set of entourages that satisfy the following property

$\forall\varepsilon>0$, $\exists\eta>0$,
$\exists N_\varepsilon\in\N$, $E^{\pm 1}$ is covered by
at most $N_\varepsilon$ sets $E_i$ such that for all $x\in M$, $E_i\circ \tilde{B}(x,\eta)$ is contained in a ball of radius $\varepsilon$.
\end{defi}

For instance, if $M$ is endowed with the discrete distance, then
$E\in \eE'$ if and only if $\forall x\in M$,
$\# E^x + \# E_x \le C$ for some $C\in \N$.

\begin{defi}
Let $M$ be a metric space. We say that $M$ satisfies property $(BG)_R$ if $\forall\varepsilon>0$, $\exists C\ge 0$ such that $\forall x\in M$, $\tilde{B}(x,R)$ is covered by at most $C$ balls of radius $\varepsilon$.
\end{defi}

We want to examine the relationship between property $(BG)_R$ and the fact that $\Delta_R\in\eE'$.

\begin{lem}\label{lem:Y1YN}
Suppose that $M$ has property $(BG)_R$. Then for all $\varepsilon>0$, there exist finitely many strictly $R$-separated subspaces $Y_1,\ldots,Y_N$ whose union $Y=Y_1\cup\cdots\cup Y_N$ is $\varepsilon$-dense.
\end{lem}

\pf
Choose $Y\subset M$ a maximal $\varepsilon$-separated subspace. By maximality, $Y$ is $\varepsilon$-dense. By property $(BG)_R$, there exists $N$ such that every ball $\tilde{B}(a,R)$ of radius $R$ is covered by $N$ balls of radius $\varepsilon/3$. Since each of these balls can contain at most one element of $Y$, $\tilde{B}(a,R)\cap Y$ has at most $N$ elements. The conclusion follows from Proposition~\ref{prop:N-R-separated}. 
\pfend

\begin{lem}\label{lem:Y1YNDelta}
Suppose that for all $\varepsilon>0$ there exists a finite union of $R$-separated subspaces $Y_1\cup\cdots\cup Y_N$ which is $\varepsilon$-dense. Then for all $R'<R/2$, $\Delta_{R'}\in \eE'$.
\end{lem}

\pf
Choose $\varepsilon>0$ and $\eta>0$ such that $2(R'+\varepsilon+\eta)<R$. Let $Y_1,\ldots,Y_N$ as in the statement of the lemma. Let $A_i=\{(y,x)\in \Delta_{R'}\vert\; \exists \tilde{y}\in Y_i,\; d(y,\tilde{y})\le\varepsilon\}$. If $a\in M$ and $(y,x),(y',x')\in A_i\circ\tilde{B}(a,\eta)$, then $d(y,y')\le d(y,x)+d(x,x')+d(x',y')\le 2R'+2\eta$, so $d(\tilde{y},\tilde{y}')<R$. Since $Y_i$ is $R$-separated, we get $\tilde{y}=\tilde{y}'$, so $d(y,y')\le 2\varepsilon$. We have shown that $A_i\circ \tilde{B}(a,\eta)$ is contained in a ball of radius $2\varepsilon$. If $E_{ij}=A_i\cap A_j^{-1}$, then $E_{ij}^{\pm 1}\circ \tilde{B}(a,\eta)$ is contained in a ball of radius $2\varepsilon$.
\pfend

\begin{lem}\label{lem:DeltaB}
If $\Delta_R\in \eE'$ then $M$ satisfies $(BG)_R$.
\end{lem}

\pf
Follows from the inclusion $\tilde{B}(x,R)\subset \Delta_R\circ\tilde{B}(x,\eta)$.
\pfend

To summarize,

\begin{prop}\label{prop:G'}
Let $M$ be a metric space. The following assertions are equivalent:
\begin{itemize}
\item[(i)] there exists $R>0$ such that $\Delta_R\in \eE'$;
\item[(ii)] there exists $R>0$ such that $M$ has $(BG)_R$;
\item[(iii)] there exists $R>0$ such that for all $\varepsilon>0$, there exists an $\varepsilon$-dense subspace $X$ such that $X$ is a finite union of $R$-separated spaces.
\end{itemize}

Moreover, if $M$ is locally compact and proper then this is equivalent to

(iv) there exists $R>0$ such that $\beta_uM$ is the union of $\bar{X}$, where $X$ runs over $R$-separated subspaces.
\end{prop}

A space that satisfies the above properties will be said to be locally of bounded geometry (LBG).

\pf
(i)$\implies$(ii): see Lemma~\ref{lem:DeltaB}.

(ii)$\implies$(iii): see Lemma~\ref{lem:Y1YN}

(iii)$\implies$(i): see Lemma~\ref{lem:Y1YNDelta}

(iv)$\implies$(iii): analogue to Proposition~\ref{prop:bounded-geometry-ULF}, (ii)$\implies$(i). Suppose that (iii) does not hold for some $\varepsilon>0$. Given any finite union of $R$-separated subspaces $X$, let $f_X(y)=(\varepsilon-d(y,X))_+$. If $X\subset Y$ then $f^{-1}_X(0)\supset f^{-1}_Y(0)$. Moreover, since $X$ is not $\varepsilon$-dense, $f^{-1}_X(0)\ne\emptyset$. By compactness, there exists $\alpha\in\beta_uM$ such that $f_X(\alpha)=0$ for all such $X$. Since $f_X=\varepsilon$ on $X$, by continuity we have $\alpha\notin\bar{X}$ (otherwise $f_X(\alpha)$ would be equal to $\alpha$). This is a contradiction.

(i)$\implies$(iv): analogue to Proposition~\ref{prop:bounded-geometry-ULF}, (i)$\implies$(ii). 
Let $\alpha\in\beta_uM$. Choose a maximal $R$-separated subspace $X$.  For all $\varepsilon>0$, there is a decomposition $\Delta_R=\cup_{i=1}^{N_\varepsilon} A_i^\varepsilon$ such that $(A_i^\varepsilon)^{\pm 1}\circ\tilde{B}(a,\eta)$ is contained in a ball of radius $\varepsilon$ for all $a\in M$. Since $\Delta_R\circ X=M$, there exists $i$ such that $\alpha\in\overline{A_i^\varepsilon\circ X}$.
Taking $\varepsilon=R/2$, there exists a family $(y_x)_{x\in X}$ satisfying $y_x\in\tilde{B}(x,R)$ such that $\alpha\in\overline{\cup_{x\in X}\tilde{B}(y_x,R/2)\cap \tilde{B}(x,R)}$. Similarly, there exist $y'_x$ such that $\alpha\in\overline{\cup_{x\in X}\tilde{B}(y'_x,R/4)\cap \tilde{B}(y_x,R/2)\cap \tilde{B}(x,R)}$, etc. We may arrange that for all $x$ and $i$, the set $Y_{i,x}=\tilde{B}(y^{(i)}_x,2^{i-1}R)\cap\cdots\cap \tilde{B}(x,R)$ is nonempty. Since $M$ is complete, there exists $z_x$ such that $\cap_i Y_{i,x}=\{z_x\}$. Let $Z=\{z_x\vert\ x\in X\}$. For all $\varepsilon>0$, $\alpha\in\overline{B(Z,\varepsilon)}$. If $\alpha\notin \bar{Z}$ then there exists a uniformly continuous function $f$ such that $f(\alpha)=1$ and $f=0$ on $Z$. By uniform continuity of $f$, there exists $\varepsilon>0$ such that $f\le 1/2$ on $B(Z,\varepsilon)$. Since $\alpha\in\overline{B(Z,\varepsilon)}$, we have $f(\alpha)\le 1/2$. Contradiction.
\pfend

In the sequel, we assume that the above properties hold. For instance, if $M$ is discrete and $\delta$-separated then
$\Delta_r\in \eE'$ for all $r<\delta$.

We remark that $\eE'$ is a coarse structure which is compatible with
the uniform structure. Moreover, every $E\in \eE'$ is contained
in an open and controlled set (for instance $\Delta_r\circ E\circ
\Delta_r$).

Let $G'(M)=\cup_{E\in \eE'}\bar{E}$.

The same proof as in Proposition~\ref{prop:G'} shows that $\exists R>0$, $\forall n\in\N$, $G(M)=\cup_X G(X)^{(n)}$, where $X$ runs over $R$-separated subspaces.

Before we prove the next proposition, we need a few lemmas.

\begin{lem}\label{lem:f-X}
Let $M$ be a locally compact metric space. Let $X\subset M$ be a closed subspace. Let $f_X=\inf(d(X,\cdot),1)$. Then $\bar{X}=f_X^{-1}(0)$ in $\beta_uM$.
\end{lem}

\pf
$\subset$ is clear. Conversely, if $\alpha\notin\bar{X}$, let us show that $f_X(\alpha)\ne 0$. There exists $f\in UC_b(M)$ such that $f_{|X}=0$ and $f(\alpha)=1$. We have $\alpha\in\overline{\{x\in M\vert\; f(x)\ge 1/2\}}$. By uniform continuity, there exists $\eta>0$ such that $d(x,X)\le\eta\implies f(x)< 1/2$. It follows that $\alpha\in
\overline{x\in M\vert\; f_X(x)>\eta\}}$, hence $f_X(\alpha)\ge\eta>0$.
\pfend

\begin{lem}\label{lem:f-alpha}
Let $M$ be a locally compact metric space. Suppose that $X,Y\subset M$ are closed subsets such that $\forall r>0$, $\exists r'>0$, $\tilde{B}(X,r')\cap\tilde{B}(Y,r')\subset \tilde{B}(X\cap Y,r)$. Then $\bar{X}\cap\bar{Y}=\overline{X\cap Y}$ in $\beta_uM$.
\end{lem}

\pf
$\supset$ is clear. Conversely, let $\alpha\in  \bar{X}\cap\bar{Y}$. Let $r>0$. Choose $r'$ as in the statement of the lemma. Since $f_X(\alpha)<r'$ and $f_Y(\alpha)<r'$, we have $\alpha\in\overline{\{x\in M\vert\; f_X(x)<r'\mbox{ and }f_Y(x)<r'\}}\subset \overline{\tilde{B}(X\cap Y,r)}\subset f_{X\cap Y}^{-1}([0,r])$. It follows that $f_{X\cap Y}(\alpha)=0$. 
\pfend

\begin{lem}\label{lem:X-cap-Y}
Let $M$ be a locally compact metric space. Let $X,Y\subset M$ be closed subsets. Then $\displaystyle\bar{X}\cap\bar{Y}=\cap_{r>0}\overline{X\cap\tilde{B}(Y,r)}$.
\end{lem}

\pf
$\subset$: let $\alpha\in\bar{X}\cap\bar{Y}$. Then
$\alpha\in\bar{X}= \overline{X\cap\tilde{B}(Y,r)\cup X\cap \tilde{B}(Y,r)^c} \overline{X\cap\tilde{B}(Y,r)}\cup \overline{X\cap \tilde{B}(Y,r)^c}$. If $\alpha$ belonged to $\overline{X\cap \tilde{B}(Y,r)^c}$, then $\alpha\in \overline{\tilde{B}(Y,r)^c}$, so $f_Y(\alpha)\ge r$. Impossible. We deduce that $\alpha\in \overline{X\cap\tilde{B}(Y,r)}$ for all $r>0$.

$\supset$: suppose $\alpha$ belongs to the right-hand side. Obviously, $\alpha\in \bar{X}$. Moreover, since $\alpha\in\overline{\tilde{B}(Y,r)}$, we have $f_Y(\alpha)\le r$ $\forall r>0$, so $f_Y(\alpha)=0$. From Lemma~\ref{lem:f-alpha}, $\alpha\in\bar{Y}$.
\pfend

\begin{prop}
Let $X$ be a closed and $\delta$-separated subset of $M$. Then $\overline{X\times X}\cap G'(M)=G'(X)\subset \beta_u(M\times M)$.
\end{prop}

\pf
$\supset$ is clear. Let us show $\subset$. We choose $\varepsilon<\delta/2$. If $g\in \overline{X\times X}\cap G'(M)$, then there exists a controlled set $A\subset M\times M$ and $\eta>0$ such that the image by $A^{\pm 1}$ of any ball of radius $\eta$ is contained in a ball of radius $\varepsilon$, and $g\in\bar{A}$. Using Lemma~\ref{lem:X-cap-Y}, for all $\varepsilon'>0$ we have $g\in\bar{B}$ where $B=(X\times X)\cap\tilde{B}(A,\varepsilon')$. We choose $\varepsilon'< \mbox{min}(\eta/2, \delta/2-\varepsilon)$. If $(x,y),(x',y)\in B$ then there exist $(a_1,a_2),(a'_1,a'_2)\in A$ such that $d(a_1,x),d(a_2,y),d(a'_1,x'),d(a'_2,y)\le\varepsilon'$. We have $d(a_2,a'_2)\le 2\varepsilon'\le\eta$, so $d(a_1,a'_1)\le 2\varepsilon<\delta$. It follows that $a_1=a'_1$, so $d(x,x')\le 2\varepsilon'+2\varepsilon<\delta$. Since $X$ is $\delta$-separated, we get $x=x'$, so the range map $r:B\to X$ is injective. Similarly, the source map $s:B\to X, (x,y)\mapsto y$ is injective. We deduce that $g\in G'(X)$.
\pfend

\begin{prop}
Let $M$ be a LBG proper metric space. Then $G'(M)$ is open in $G(M)$, thus is locally compact. Moreover, it has a Haar system.
\end{prop}

\pf
Let $E\in\eE'$. Let $r>0$ such that $\Delta_r\in \eE'$, and let $E'=\Delta_r\circ E\circ \Delta_r$. It suffices to prove that $\overline{E'}$ is a neighborhood of $\bar{E}$ in $G(M)$. This follows from
$\bar{E}\subset f_E^{-1}([0,r/3])\subset f_E^{-1}([0,r/2))\subset\overline{E'}$ (see notation in Lemma~\ref{lem:f-X}).

The proof of the last assertion is almost the same as in the case of a ULF space, so we omit it.
\pfend

The drawback of the groupoid $G'(M)$ is that if $X\subset M$ is $R$-dense then the inclusion $G'(X)\to G'(M)$ is not necessarily a Morita equivalence. To remedy this, we define

\begin{defi}\label{def:G(M)}
Let $M$ be a LBG, locally compact proper metric space. We define $G(M)$ as the union of all $\bar{E}$, where $E\in\eE'$ and $r(E)$, $s(E)$ have bounded geometry.
\end{defi}

An alternative definition is : $G(M)=\cup_{X}G(X)$, where $X$ runs over closed, BG subspaces.

\begin{lem}\label{lem:E-X}
Let $M$ be a metric space. If $X\subset M$ has bounded geometry and $E\in \eE'$, then $E_X$ and $E^X$ have bounded geometry.
\end{lem}

\pf
We prove the first assertion, the second being similar.
Let $R>0$ and $\varepsilon>0$. We want to show that there exists $n$ such that every ball (in $X$) of radius $R$ can be covered by $n$ balls of radius $\varepsilon$. Let $R'$ such that $E\subset \Delta_{R'}$. Let $\eta>0$ such that $\exists N$, $\forall a\in M$, $E\circ \tilde{B}(a,\eta)$ can be covered by $N$ balls of radius $\varepsilon$.

Let $y\in E_X$. There exists $x\in X$ such that $(y,x)\in E$. For all $y'\in Y\cap \tilde{B}(y,R)$, there exists $x'\in X$ such that $(y',x')\in E$. Then $d(x,x')\le R+2R'$. Now, there exists $N'$ (dependent on $\eta$ and $R+2R'$) such that $X\cap \tilde{B}(x,R+2R')$ can be covered by $N'$ balls $B_i$ (on $X$) of radius $\eta$. Since $y'\in\cup_i E\circ B_i$, $\tilde{B}(y,R)$ can be covered by $NN'$ balls of radius $\varepsilon$.
\pfend

Let us denote $\beta'_uM=\cup\bar{X}$, where $X$ runs over all bounded geometry subspaces $X$.

\begin{prop}
Let $M$ be a LBG, proper metric space. Then $\beta'_uM$ is an open subspace of $\beta_uM$ which is saturated for the action of $G'(M)$.
\end{prop}

\pf
Let $r>0$ such that $\Delta_r\in \eE'$. For all $E\in\eE'$ such that $s(E)$ and $r(E)$ have bounded geometry, $E'=\Delta_r\circ E\circ \Delta_r$ belongs to $\eE'$ and $s(E')$, $r(E')$ have bounded geometry thanks to Lemma~\ref{lem:E-X}. Therefore, $\overline{E'}$ is a neighborhood of $E$ in $G'(M)$. We deduce that $\beta'_uM$ is open.

Let us show that $\beta'_uM$ is saturated. Let $g\in G'(M)$ such that $s(g)\in\beta'_uM$. We have to show that $r(g)\in\beta'_uM$.

There exists $E\in\eE'$ such that $g\in\bar{E}$. Moreover, there exists a bounded geometry subspace $X$ such that $s(g)\in\bar{X}$.

By Lemma~\ref{lem:X-cap-Y}, $s(g)=\overline{s(E)}\cap\bar{X}
=\cap_{r'>0}\overline{s(E)\cap\tilde{B}(X,r')}$. Since $g\in\bar{E}
=\overline{E\cap s^{-1}(\tilde{B}(X,r'))}\cup\overline{E\cap s^{-1}(\tilde{B}(X,r)^c)}$, we must have $g\in\bar{E}
=\overline{E\cap s^{-1}(\tilde{B}(X,r'))}$ (otherwise $s(g)\in \overline{\tilde{B}(X,r'/2)} \cap \overline{\tilde{B}(X,r')^c}=\emptyset$ (see Lemma~\ref{lem:X-cap-Y}).

After replacing $E$ by $E\cap s^{-1}(\tilde{B}(X,r))$, we may assume that $s(E)$ has bounded geometry (since $\tilde{B}(X,r)=\Delta_r\circ X$), so $r(E)$ also has bounded geometry (see Lemma~\ref{lem:E-X}). We deduce that $r(g)\in r(E)\subset \beta'_uM$.
\pfend

From this, we deduce easily

\begin{prop}
Let $M$ be a proper, LBG metric space. Then $G(M)=G'(M)_{\beta'_uM}$ is a locally compact groupoid with Haar system.
\end{prop}

Remark : $G(M)$ is generally not $\sigma$-compact if $M$ does not have bounded geometry.

\begin{prop}\label{prop:subspace-Morita}
Let $M$ be a proper, LBG metric space. If $r>0$ is such that $\Delta_{2r}\in\eE'$, then given any closed, $r$-dense subspace $N$, the inclusion $G(N)\to G(M)$ is a Morita equivalence.
\end{prop}

\pf
Indeed, $\beta'_uN$ is a closed transversal for $G(M)$. Since $G(N)=G(M)_{\beta'_uN}^{\beta'_uN}$, we get the result.
\pfend

\section{The classifying space for proper actions of an \'etale groupoid}

In this section, $G$ denotes a locally compact, $\sigma$-compact, \'etale groupoid. Given a compact subset $K$ of $G$, let $P_K(G)$ be the space of probability measures $\mu$ on $G$ such that for all $g,h\in\mbox{supp}(\mu)$, $r(g)=r(h)$ and $g^{-1}h\in K$. We endow $P_K(G)$ with the weak-* topology, and the natural left action of $G$.
Note that the support of $\mu$ must be finite, as it is discrete and included in a compact set of the form $C(g)= \{gk\vert\; k\in K,\;r(k)=s(g)\}$.

\begin{prop}
The action of $G$ on $P_K(G)$ is proper and cocompact.
\end{prop}

\pf
Let us show that the action is proper. If $L$ is a compact subset of $G$, it is a standard exercise to check that the set $C_L=\{\mu\in P_K(G)\vert\;\mbox{supp}(\mu)\subset L\}$ is an exhausting sequence of compact subsets of $G$. Now, if $\mu\in C_L$ and $g\mu\in C_L$, then $g$ belongs to the compact set $LL^{-1}=\{hk^{-1}\vert\ h,k\in L\}$, so the action is proper.

The action is cocompact since the saturation of $C_K$ is equal to $P_K(G)$.
\pfend

\begin{lem}\label{lem:PK=EG}
Let $Y$ be a proper and $G$-compact $G$-space. Then there exists a compact subset $K$ of $G$ and a continuous equivariant map $Y\to P_K(G)$.
\end{lem}

\pf
Since the action of $G$ on $Y$ is proper, there exists $c\in C_c(Y)_+$ such that $\sum_g c(yg)=1$. Let $\mu_y=\sum_g c(yg)\delta_g$. Let $L$ be the support of $c$. There exists a compact subset $K$ of $G$ such that $\forall (y,g)\in Y\times_{G^{(0)}}G$, $(y,yg)\in L\times L\implies g\in K$. Then for all $g,h\in\mbox{supp}(\mu_y)$, we have $g^{-1}h\in K$, so $y\mapsto \mu_y$ determines an equivariant map $Y\to P_K(G)$.
\pfend

Before we proceed, we need a few lemmas.

\begin{lem}\label{lem:h(b)}
Let $a,a',b$ be selfadjoint elements of an abelian $C^*$-algebra, and $\varepsilon>0$. Suppose $a'(1-a)=0$, $-1\le b\le 1$ and $\|a(1-b^2)\|\le\varepsilon$.

Let $h:[-1,1]\to [-1,1]$ continuous such that $h(0)=0$, $h(t)=-1$ on $[-1,-1+\sqrt{1-\varepsilon}]$, and $h(t)=1$ on $[1-\sqrt{1-\varepsilon},1]$. Let $b'=h(b)$. Then $a'(1-{b'}^2)=0$.
\end{lem}

\pf
We may assume that the $C^*$-algebra is $C(X)$, where $X$ is a compact space. After evaluating at each point, we may assume that $a,a',b$ are real numbers. If $a'\ne 0$ then $a=1$, so $|1-b^2|\le\varepsilon$, so $b'=\pm 1$.
\pfend

\begin{lem}\label{lem:ABJ}
Let $A$ and $B$ be $G$-algebras, $J$ a $G$-invariant ideal of $B$. Suppose that  $[(\eE,F)]\in KK_G(B,A)$ satisfies
\begin{itemize}
\item[(i)] $j(F^2-1)=0$ for all $j\in J$
\item[(ii)] $[b,F]=0$ for all $b\in B$.
\end{itemize}
Let $\eE'=\{x\in \eE\vert\; Jx=0\}$. Then $F$ induces $F'\in\lL(\eE')$, and $(\eE',F')$ determines an element of $KK_G(B/J,A)$ whose image in $KK_G(B,A)$ is equal to $[(\eE,F)]$.
\end{lem}

\pf
The first assertion comes from the fact that $F$ commutes with $J$.

Since $BJ\subset J$, $B$ maps to $\lL(\eE')$, and this maps obviously factors through $B/J$.

Clearly, $B$ commutes with $F'$. It remains to check that $B({F'}^2-1)$ is compact. Let $T={F}^2-1$. Let $b=b^*\in B$.

$(bT)^3=T(b^3T)T\in T\kK(\eE) T=\overline{\mbox{span}}\{\theta_{T\xi,T\eta}\vert\;\xi,\eta\in\eE\}$ (where $\theta_{\xi,\eta}$ denotes the rank-one operator $\zeta\mapsto \xi\langle\eta,\zeta\rangle$). Now, $T\xi$, $T\eta\in\eE'$, so $(bT)^3$ induces an element of $\kK(\eE')$. Taking the cubic root, we get that $b({F'}^2-1)$ is compact.
\pfend

\begin{defi}
A map $f:X\to Y$ between two topological spaces is said to be locally injective if $X$ is covered by open subsets $U$ for which $f_{|U}$ is injective.
\end{defi}

If $Z$ is a proper $G$-space and $A$ is a $G$-algebra, we denote by $RK_G(Z;A)$ the inductive limit of $KK_G(C_0(Y),A)$, where $Y$ runs over $G$-compact subspaces of $Z$.

\begin{lem}\label{lem:loc-injective}
Let $G$ be a locally compact \'etale groupoid.
Let $Y$ and $T$ be locally compact spaces endowed with an action of $G$, such that the action of $G$ on $Y$ is proper and cocompact. Assume that the map $p:Y\to G^{(0)}$ is locally injective. Then the natural map $RK_{T\rtimes G}(T\times Y;A)\to RK_G(Y;A)$ induced by the second projection $T\times Y\to Y$ is an isomorphism.
\end{lem}

\pf
We want to construct a map in the other direction. Let $[(\eE,\varphi,F)]$ be an element of $RK_G(Y;A)$.

Let $K$ a compact subset of $Y$ such that $KG=Y$. There exists a finite open cover $(U_i)$ of $K$ for which $p_{|U_i}$ is injective. There exist $f_i\in C_c(Y)_+$ such that $\mbox{supp}(f_i)\subset U_i$ and $K\subset f_i^{-1}((0,+\infty))$. After replacing $f_i(y)$ by $f_i(y)/\sum_{j,g}f_j(yg)$, we can assume that $\sum_{j,g}f_j(yg)=1$ for all $y\in Y$.

Consider $F'_x=\sum_{i,g\in G^x} \alpha_g(f_i^{1/2}F_{s(g)}f_i^{1/2})$. By construction, $F'$ is a self-adjoint and $G$-invariant operator. Let us check that it is a compact perturbation of $F$.

\begin{eqnarray*}
h(F'_x-F_x) &=&
\sum_{i,g} g\left(
\alpha_g(f_i^{1/2}F_{s(g)}f_i^{1/2}) -\alpha_g(f_i^{1/2}\alpha_g(f_i^{1/2})\alpha_g(F_{s(g)})\right)\\
&&
+h\alpha_g(f_i)(\alpha_g(F_{s(g)})-F_{r(g)}).
\end{eqnarray*}

Let $L=\{g\in G\vert\;\exists i,\;\exists y\in\mbox{supp}(h),\; f_i(yg)\ne 0\}$. Then $L$ is relatively compact, and the term in the sum is zero when $g\notin L$, so for each $x$ the sum is finite. In addition, each term is compact, so the sum is compact.

By local injectivity of $Y\to G^{(0)}$, $f_i^{1/2}F_{s(g)}f_i^{1/2}$ commutes with $C_0(Y_{s(g)})$, so $\alpha_g(f_i^{1/2}F_{s(g)}f_i^{1/2})$ commutes with $C_0(Y_{r(g)})$ (where $Y_x$ denotes the fiber of $Y$ over $x\in G^{(0)}$). Therefore, $F'$ commutes with $C_0(Y)$.

After replacing $F$ by $F'$, we can assume that $F$ is $G$-invariant and commutes with $C_0(Y)$. Since $A$ is a $C_0(T)$-algebra, $F$ also commutes with the action of $C_0(T)$, so $F$ is an endomorphism of the left $C_0(T\times_{G^{(0)}} Y)$-module $\eE$. We can also assume that $-1\le F\le 1$.

Let $f,f'\in C_c(Y)_+$ such that $f'=1$ on $K$ and $f=1$ on the support of $f'$. Let $\varepsilon\in (0,1)$. Since $f(1-F^2)$ is compact, there exists a compact set $L\subset T$ such that $\|f(1-F^2)_t\|\le \varepsilon$ for all $t\in T-L$. Let $F'=h(F)$ where $h$ is like in Lemma~\ref{lem:h(b)}, then $f'(1-{F'}^2)_t=0$ for all $t\notin L$.

Let us show that $\varphi(1-F^2)=0$ for all $\varphi\in C_c(T\times_{G^{(0)}}Y)$ supported outside $L\times_{G^{(0)}}K$. Since $\varphi$ is a finite sum of functions $\varphi_i$ supported in sets of the form $U\times_{G^{(0)}}V$, where $U$ and $V$ are open, relatively compact sets, which are domains of local homeomorphisms coming from some element $g_i\in G$ such that $(U\times_{G^{(0)}} V)g\subset L^c\times \mbox{supp}(f')$, we may assume that $\varphi$ is equal to one of those $\varphi_i$'s. Choose $h_1\in C_c(T)_+$ and $h_2\in C_c(Y)_+$ such that $U=h_1^{-1}(\R_+^*)$ and $V=h_2^{-1}(\R_+^*)$. Since $(t,y)\mapsto h_1(tg^{-1})h_2(yg^{-1})$ is zero outside $L^c\times\mbox{supp}(f')$, we have $g\cdot(h_1\otimes h_2)(1-{F'}^2)=0$. By $G$-invariance of $F'$, we have $(h_1\otimes h_2)(1-{F'}^2)=0$. We deduce that $\varphi(1-{F'}^2)=0$.

Now, let $Y'$ be the saturation of $L\times_{G^{(0)}}K$. Using Lemma~\ref{lem:ABJ} for $B=C_0(T\times_{G^{(0)}}Y)$ and $J=C_0(T\times_{G^{(0)}}Y-Y')$, we get an element of $RK_G(Y';A)$. In fact, the construction of Lemma~\ref{lem:ABJ} yields an element of $RK_{T\rtimes G}(Y';A)$, and one easily checks that the map $RK_G(Y;A)\to RK_{T\rtimes G}(Y';A)$ is inverse to the map $RK_{T\rtimes G}(Y';A) \to RK_G(Y;A)$.
\pfend

\begin{defi}
Let $G$ be a locally compact groupoid. A $G$-simplicial complex of dimension $\le n$ is a pair $(X,\Delta)$ given by
\begin{itemize}
\item[(i)] a locally compact space $X$ (the set of vertices), with an action of $G$ relative to a locally injective map $p:X\to G^{(0)}$;
\item[(ii)] a closed, $G$-invariant subset $\Delta$ of the space of measures on $X$ (endowed with the weak-$*$ topology), such that each element of $\Delta$ is a probability measure whose support (called a simplex) has at most $n+1$ elements and is a subset of one of the fibers of $p$. In addition, we require that if $\mbox{supp}(\mu)\subset\mbox{supp}(\nu)$ and $\nu\in \Delta$, then $\mu\in\Delta$.
\end{itemize}

The $G$-simplicial complex is typed if there is a discrete set $\tT$ (the set of types) and a $G$-invariant, continuous map $\tau:X \to\tT$ such that the restriction of $\tau$ to any simplex is injective.
\end{defi}

It is not hard to see that $\Delta$ is locally compact, and that if $G$ acts properly on $X$ then it acts properly on $\Delta$.

The barycentric subdivision $(X',\Delta')$ is the $G$-simplicial complex whose vertex set consists of the centers of simplices of $\Delta$, such that $S=\{\nu_0,\ldots,\nu_k\} $ is a simplex if and only if the union of the supports of $\nu_i$ is a simplex of $\Delta$. Using local injectivity of $p$, we see that $X'$ is a closed subspace of $\Delta$, so that $G$ acts properly on $X'$. It is clear that $X'\to G^{(0)}$ is also locally injective.

This construction shows that if a $G$-space has a structure of $G$-simplicial complex, then it has the structure of typed $G$-simplicial complex.

Let us introduce the following notation: if $A$ is a $G$-algebra, then $\mbox{BC}(G;A)$ means that $G$ satisfies the Baum--Connes conjecture with coefficients in $A$.

We now prove the following generalization of \cite{ceo}:

\begin{theo}\label{thm:BCT}
Let $G$ be a locally compact, second-countable \'etale groupoid, $T$ a locally compact, second-countable $G$-space, and $A$ a $T\rtimes G$-algebra. Then the canonical map
$$K_*^{top}(T\rtimes G;A)\to K_*^{top}(G;A)$$
is an isomorphism. As a consequence, $\mbox{BC}(G;A) \iff \mbox{BC}(T\rtimes G;A)$
\end{theo}

\pf
This amounts to showing that $RK_{T\rtimes G}(T\times_{G^(0)}Z;A)\to RK_G(Z;A)$ is an isomorphism when $Z=\underline{E}G$ is the classifying space for proper actions of $G$. Using Lemma~\ref{lem:PK=EG}, it suffices to prove this for $Z=P_K(G)$. Since $P_K(G)$ is a proper, $G$-compact $G$-simplicial complex of dimension $n=\sup\{\#G^x_K\vert\;x\in K\}$, it suffices to show the isomorphism for any typed, proper $G$-compact $G$-simplicial complex.

We proceed by induction on the dimension $n$. For $n=0$, this is the content of Lemma~\ref{lem:loc-injective}.

Suppose the result is true in dimensions $<n$. Let $Z$ by a typed, proper $G$-compact $G$-simplicial complex of dimension $n$ and let $F$ be its $n-1$-skeleton. Let $U=Z-F$. Consider the diagram

$$\xymatrix{
\cdots\ar[d] & \cdots\ar[d] \\
RK_{T\rtimes G}(T\times_{G^{(0)}}F,A)\ar[r]\ar[d] &
RK_G(F,A)\ar[d]\\
RK_{T\rtimes G}(T\times_{G^{(0)}} Z,A)\ar[r]\ar[d]  &
RK_G(Z,A)\ar[d]\\
\lim_Y  KK_G(C_0(Y\cap (T\times_{G^{(0)}}U),A)\ar[r]\ar[d]  &
KK_G(C_0(U),A)\ar[d]\\
\cdots & \cdots
}
$$
where the inductive limit is over $G$-compact subspaces of $T\times_{G^{(0)}} Z$.

The columns are exact thanks to Lemma~\ref{lem:K-exact} below.

The first horizontal arrow is an isomorphism, thanks to the induction assumption, so in order to use the five-lemma, we need to show that
$$ \lim_Y  KK_G(C_0(Y\cap (T\times_{G^{(0)}}U)),A)\to KK_G(C_0(U),A)$$
is an isomorphism.

Let $Z'$ be the set of centers of $n$-simplices. Since the simplicial complex is typed, $U$ is isomorphic to $Z'\times \R^n$, so the right-hand side is $RK_G^n(Z';A)$.

Intersections of $T\times_{G^{(0)}}U$ by $G$-compact subspaces of $T\times_{G^{(0)}}Z$ correspond to intersections of $T\times_{G^{(0)}} Z' \times \R^n$ with closed subsets of sets of the form $Y'\times \R^n$ with $Y'\subset T\times_{G^{(0)}} Z'$ $G$-compact. Therefore, the left-hand side is isomorphic to $\lim_{Y'} KK_{T\rtimes G}(C_0(Y'\times \R^n),A)=RK_{T\rtimes G}^n(T\times_{G^{(0)}}Z';A)$. The assertion follows from Lemma~\ref{lem:loc-injective}.
\pfend

In the proof of the above theorem, we used :

\begin{lem}\label{lem:K-exact}
Let  $G$ is a locally compact, second-countable groupoid with Haar system, $Z$ is a second-countable, proper $G$-space, $F$ a $G$-invariant subset of $Z$ and $U$ its complementary, then for any $G$-algebra $A$ there is a six-term exact sequence
$$\xymatrix{
KK_G(C_0(F),A)\ar[r] & KK_G(C_0(Z),A)\ar[r] & KK_G(C_0(U),A)\ar[d] \\
KK_G^1(C_0(U),A)\ar[u] & KK_G^1(C_0(Z),A)\ar[l] & KK_G^1(C_0(F),A)\ar[l] 
}$$
\end{lem}

\pf
This is a consequence of \cite[Corollaire~5.2]{tu99a} and the proof of \cite[Th\'eor\`eme~5.2]{m99}.
\pfend

\begin{theo}\label{thm:BCcoef1}
Let $H\le G$ be locally compact, second-countable \'etale groupoids (with $H$ closed in $G$). If $G$ satisfies the Baum-Connes conjecture with coefficients, then $H$ satisfies the Baum-Connes conjecture with coefficients.
\end{theo}

\pf
The proof is the same as in \cite{ceo}.
\pfend

Now, we want to remove the second-countability assumption.

\begin{lem}
Let $G$ be a locally compact, second-countable, \'etale groupoid. Suppose that $G$ acts on a $\sigma$-unital $C^*$-algebra $A$. Then there exists a sub-$C^*$-algebra $B$, invariant by $G$, which contains an approximate unit for $A$.
\end{lem}

\pf
Let $(U_i)$ and $(U'_i)$ be countable families of open subsets of $G$ such that $U_i\subset U'_i$, $G=\cup_i U_i$, $\{U_i\}$ is stable by the inversion map, and $r$ and $s$ induce homeomorphisms from $U'_i$ onto their respective images $W'_i$ and $V'_i$. Let $V_i=s(U_i)$ and $W_i=r(U_i)$. View $U_i$ as a homeomorphism from $V_i$ to $W_i$. We consider $A$ as an upper semi-continuous field of $C^*$-algebras over $G^{(0)}$. Denote by $A_{V_i}$ the set of restrictions of elements of $A$ to $V_i$. Then $U_i$ induces an isomorphism from $A_{V_i}$ to $A_{W_i}$.

Let $(a_n)$ be a countable approximate unit of $A$. Let $X_0=\{a_n\vert\; n\in\N\}$. Let $X_1$ be a countable subset of $A$ such that for all $i$ and all $n$, ${a_n}_{|V_i}$ is the restriction to $W_i$ of an element of $X_1$. In the same way, we define $X_2$, $X_3$, etc. Then the $C^*$-algebra $B$ generated by $\cup_n X_n$ satisfies the required properties.
\pfend

\begin{prop}\label{prop:TT'}
Let $G$ be a locally compact, second-countable, \'etale groupoid. Suppose that $G$ acts on a locally compact, $\sigma$-compact space $T$. Then there exists a locally compact, second-countable space $T'$ with an action of $G$, and a continuous proper equivariant map $T\to T'$ with dense image.
\end{prop}

\pf
We apply the preceding lemma to $A=C_0(T)$. There exists a $G$-invariant separable subalgebra $B$ containing an approximate unit. If $T'$ denotes the spectrum of the abelian $C^*$-algebra $B$, then the inclusion $B\to A$ induces a map $T\to T'$ with the required properties.
\pfend

\begin{theo}
Let $G$ be a locally compact groupoid isomorphic to $X\rtimes G'$, where $X$ is a locally compact, $\sigma$-compact space and $G'$ is a locally compact, second-countable and \'etale groupoid. We assume that the anchor map $X\to {G'}^{(0)}$ of the action is proper. Let $T$ be a locally compact, $\sigma$-compact $G$-space. Let $A$ be a $T\rtimes G$-algebra. Then $\mbox{BC}(G;A) \iff \mbox{BC}(T\rtimes G;A)$.
\end{theo}

\pf
Let $T'$ as in Proposition~\ref{prop:TT'}. We note that $\mbox{BC}(T\rtimes G';A) \iff \mbox{BC}(T'\rtimes G';A)$. To see this, we have to show that the forgetful functor from $RK_{T\rtimes G'}(Y;A)$ to $RK_{T'\rtimes G'}(Y';A)$ is an isomorphism, where $Y'\subset T'\times_{G^{(0)}} \underline{E}G$ is $G$-compact and $Y$ is its inverse image in $T\times_{G^{(0)}} \underline{E}G$. (Note that $Y$ is $G$-compact thanks to the properness assumption of $Y\to Y'$).
Let $[(E,\varphi,F)]$ be an element of $RK_{T'\rtimes G'}(Y';A)$. Let $(f_i)$ be an approximate unit of $C_0(T')$. For all $h\in C_0(\underline{E}G)$, the operator $\varphi(f_i\otimes h)(F^2-1)$ is compact. It follows that $\varphi(f\otimes f)(F^2-1)$ is compact for all $f\in C_0(T)$ (since $f_i$ is an approximate unit in $C_0(T)$). Using a similar argument for $[f_i\otimes h,F]=f_i [h,F]$, we find that $[(E,\varphi,F)]$ determines an element of $RK_{T\rtimes G'}(Y;A)$.

This proves the assertion $\mbox{BC}(T\rtimes G';A) \iff \mbox{BC}(T'\rtimes G';A)$. Similarly, $\mbox{BC}(G';A) \iff \mbox{BC}(G^{(0)}\rtimes G';A)$. Moreover, from Theorem~\ref{thm:BCT}, we get $\mbox{BC}(G';A) \iff \mbox{BC}(T'\rtimes G';A)$. Combining these three equivalences, we get the conclusion.
\pfend

As above, we deduce the following generalization of Theorem~\ref{thm:BCcoef1}.

\begin{theo}\label{thm:BCcoef}
Let $G$ be a locally compact groupoid isomorphic to $X\rtimes G'$, where $X$ is a locally compact, $\sigma$-compact space and $G'$ is a locally compact, second-countable and \'etale groupoid. We assume that the anchor map $X\to {G'}^{(0)}$ of the action is proper. Let $H$ be a closed, \'etale subgroupoid of $G$. If $G$ satisfies the Baum-Connes conjecture with coefficients, then $H$ also satisfies the Baum-Connes conjecture with coefficients.
\end{theo}

\section{The coarse Baum--Connes conjecture with coefficients}
\begin{defi}
Let $M$ be a LBG, proper metric space. We say that $M$ satisfies the coarse Baum--Connes conjecture (resp. the coarse Baum--Connes with coefficients) if the groupoid $G(M)$ satisfies the Baum--Connes conjecture with coefficients in the $C^*$-algebra $UC_b(M,\kK)$ (resp. the Baum--Connes conjecture with arbitrary coefficients). 
\end{defi}

We define analogously the full coarse Baum--Connes conjecture (see also \cite{oy09}) and the full coarse Baum--Connes conjecture with coefficients.

Note that if $r>0$, and if $X$ is a maximal $r$-separated subspace, then $\beta X$ is a complete transversal, so $G(X)$ is Morita equivalent to $G(M)$. Since $\mbox{BC}(G(M);UC_b(M,\kK))$ is equivalent to $\mbox{BC}(G(X);\ell^\infty(X,\kK))$, the coarse Baum--Connes conjecture for $M$ coincides with the usual one \cite{sty} when $M$ has bounded geometry.

\begin{theo}
Let $M$ be a bounded geometry, proper metric space and $N$ a closed subspace. If $M$ satisfies the coarse Baum--Connes conjecture with coefficients then $N$ also does. A similar assertion holds for the full coarse Baum--Connes conjecture.
\end{theo}

\pf
Let $X\subset M$ be a maximal $1$-separated subspace. Let $Y=\{x\in X\vert\; d(x,N)\le 1\}$. Since $G(X)$ is Morita equivalent to $G(M)$, it satisfies $\mbox{BC}_{coef}$. From \cite{sty}, the groupoid $G(X)$ satisfies the conditions of Theorem~\ref{thm:BCcoef}. Since $G(Y)$ is a closed subgroupoid of $G(X)$, it also satisfies $\mbox{BC}_{coef}$. Finally, $G(N)$ satisfies $\mbox{BC}_{coef}$ since it is Morita equivalent to $G(Y)$.

The proof for the full coarse Baum--Connes conjecture is similar.
\pfend

Our goal is now to examine the question of finding a ``descent principle''. It is known that if $\Gamma$ is a (torsion free) discrete group whose classifying space $B\Gamma$ is a finite CW-complex, the coarse Baum--Connes conjecture for the underlying metric space of $\Gamma$ implies that the Baum--Connes map for the group $\Gamma$ is injective, but it is not known whether one can extend this descent principle to more general groups (such as groups with torsion such that $\underline{E}\Gamma$ is $\Gamma$-compact). One might wonder whether coarse Baum--Connes conjecture with coefficients is strong enough to imply injectivity of the Baum--Connes map for the group. Since we are not able to answer this question, we introduce the following definition:

\begin{defi}
Let $M$ be a LBG, proper metric space. We say that $M$ satisfies the strong coarse Baum--Connes conjecture with coefficients (SCBC) if for all $n\in\N^*$, the semi-direct product groupoid $G(M)^n\rtimes S_n$ of $G(M)^n$ by the symmetric group $S_n$ (acting by permutation on the factors of $G(M)^n$) satisfies the Baum--Connes conjecture with coefficients. One defines analogously the strong full coarse Baum--Connes conjecture with coefficients (SFCBC).
\end{defi}

For instance, if $M$ admits a coarse embedding into a Hilbert space, then $M$ satisfies (SCBC) and (SFCBC). Indeed, we can reduce to the case when $M$ is discrete. Then $G(M)$ acts properly on a continuous field of Hilbert spaces. It follows immediately that $G(M)^n\rtimes S_n$ also does, so that it satisfies $\mbox{BC}_{coef}$ by \cite{tu99a} (see \cite{sty} for $n=1$).

Before we state the next theorem, we note that if $F$ is a finite group acting by isometries on a LBG space $M$, there is an obvious notion of $F$-equivariant coarse Baum--Connes conjecture (the coarse assembly map taking its values in $K_F(C^*(M))\cong K(C^*(M)\rtimes F)$), which is shown (by essentially the same methods) to be equivalent to the Baum--Connes conjecture for $G(M)\rtimes F$ with coefficients in $UC_b(M,\kK)$. When $M$ is the underlying metric space of a discrete group endowed with any left-invariant proper distance, this is again equivalent to the Baum--Connes conjecture for $\Gamma$ with coefficients in $\ell^\infty(\Gamma,\kK)\rtimes F$, where the actions of $\Gamma$ and of $F$ on $\ell^\infty(\Gamma,\kK)$ are induced by the right action of $\Gamma$ and the left action of $F$ on $\Gamma$.

\begin{theo}
Let $\Gamma$ be a countable group, and let $X$ be the underlying metric space (given by any left-invariant proper distance). Consider the following statements:
\begin{itemize}
\item[(i)] $X$ satisfies SFCBC;
\item[(ii)] for every finite subgroup $F$ of $\Gamma$, $X$ satisfies the $F$-equivariant full coarse Baum--Connes conjecture;
\item[(iii)] the full Baum--Connes map for $\Gamma$ is injective.
\end{itemize}
Then $(i)\implies (ii)$. Moreover, if there is a classifying space for proper actions $\underline{E}\Gamma$ which is second-countable and $\Gamma$-compact, and if $\underline{E}\Gamma$ is $F$-equivariantly uniformly contractible for every finite subgroup $F$ of $\Gamma$, then $(ii)\implies (iii)$.
\end{theo}

(A space is said to be uniformly contractible if there is a uniformly continuous homotopy between the identity map and the constant map.)

\pf
To prove $(i)\implies (ii)$, we note that if $n=\# F$ then there are embeddings $G(X)\rtimes F\to G(X)^F\rtimes F \to G(X)^F\rtimes S_F\cong G(X)^n\rtimes S_n$, where $G(X)^F$ denotes the set of maps from $F$ to $G(X)$, the first map being given by the $F$-equivariant embedding $G(X)\to G(X)^F$, $\gamma\mapsto (f^{-1}(\gamma))_{f\in F}$. Property $(ii)$ follows from (the full version of) Theorem~\ref{thm:BCcoef}.

Let us prove $(ii)\implies (iii)$. We first prove that for every proper finite $\Gamma$-simplicial complex $Y$, the full Baum--Connes conjecture $\mbox{FBC}(\Gamma;UC_b(Y,\kK))$ holds. We recall that, up to uniform-coarse equivalence, there is one and only one distance on $Y$ which is $\Gamma$-invariant and proper.

If $Y$ is 0-dimensional, it is isomorphic to $\Gamma/F$ where $F$ is a finite group, so this reduces to $\mbox{BC}(\Gamma;\ell^\infty(\Gamma/F,\kK))$, which is true by (ii).

If $Y$ is arbitrary, we may assume that $Y$ is a typed $\Gamma$-simplicial complex. We proceed by induction on the dimension of $Y$. Let $n$ be the dimension of $Y$ and suppose the result is true in dimensions $<n$. Let $Y'$ be the $n-1$-skeleton of $Y$ and $U$ its complementary. We have the exact sequence
$$0\to UC_{b,0}(U,\kK)\to UC_b(Y,\kK)\to UC_b(Y',\kK)\to 0,$$
where $UC_{b,0}(U,\kK)$ denotes the algebra of uniformly continuous and bounded functions from $U$ to $\kK$ that vanish at infinity.
We note that $UC_{b,0}(U,\kK)\cong UC_b(Y'',\kK)\otimes C_0(\R^n)$ where $Y''$ is the set of centers of the open simplices in $U$.

Taking the full crossed-product with $\Gamma$ preserves exact sequences, so we get a six-term exact sequence
$$\cdots\to K_{i+n}(UC_b(Y'',\kK)\rtimes_f\Gamma)
\to K_i(UC_b(Y,\kK)\rtimes_f\Gamma)\to K_i(UC_b(Y',\kK)\rtimes_f\Gamma)\to\cdots$$

Similarly, let us show that we have exact sequences in topological $K$-theory:

$$\cdots\to K_{i+n}^{top}(\Gamma;UC_b(Y'',\kK))
\to K_i^{top}(\Gamma;UC_b(Y,\kK))\to K_i^{top}(\Gamma;UC_b(Y',\kK))\to\cdots$$
Indeed since $K_*^{top}(\Gamma;A)$ is the inductive limit of $KK_\Gamma(C_0(P_d(\Gamma)),A)$, we have to check that for every proper and finite $\Gamma$-simplicial complex $Z$,
$$\cdots\to KK_\Gamma^n(C_0(Z);UC_b(Y'',\kK))
\to KK_\Gamma(C_0(Z);UC_b(Y,\kK))\to KK_\Gamma(C_0(Z);UC_b(Y',\kK))\to\cdots$$
is exact. When $Z$ is 0-dimensional, it is isomorphic to $\Gamma/F$ where $F$ is a finite group, so by \cite[Proposition~5.14]{ce} this reduces to
$$\cdots\to K^F_n(UC_b(Y'',\kK))
\to K^F(UC_b(Y,\kK))\to K^F(UC_b(Y',\kK))\to\cdots$$
which is indeed exact since the functor $K^F$ preserves exact sequences.

When $Z$ is arbitrary, this follows from an induction on the dimension of $Z$ and from Lemma~\ref{lem:K-exact}.

Since the Baum--Connes assembly map intertwines the two above exact sequences, an application of the five-lemma completes the proof that $\mbox{FBC}(\Gamma;UC_b(Y,\kK))$ holds. 
Now, since $\underline{E}\Gamma$ is compact, there exists $Y$ of the form $P_d(\Gamma)$ and equivariant maps $\underline{E}\Gamma\to Y\to \underline{E}\Gamma$ (see Lemma~\ref{lem:PK=EG}) whose composition is $\Gamma$-homotopic to the identity. We thus get maps
$UC_b(\underline{E}\Gamma,\kK) \to UC_b(Y,\kK) \to UC_b(\underline{E}\Gamma,\kK)$ whose composition is $\Gamma$-homotopic to the identity.
It follows that the full Baum--Connes assembly map for $\Gamma$ with coefficients in $UC_b(\underline{E}\Gamma,\kK)$ is a direct factor of the full Baum--Connes assembly map for $\Gamma$ with coefficients in $UC_b(Y,\kK)$, so it is an isomorphism.

Now, consider the diagram
$$\xymatrix{
K^{top}(\Gamma)\ar[r]\ar[d] & K(C^*(\Gamma))\ar[d] \\
K^{top}(\Gamma;UC_b(\underline{E}\Gamma;\kK))\ar[r] & K(UC_b(\underline{E}\Gamma;\kK))
}
$$
We have shown that the lower horizontal arrow is an isomorphism.
The leftmost vertical arrow is an isomorphism. To see this, using again exact sequences and an induction argument, we are are reduced as above to showing that $K_F(\C)\to K_F(UC_b(\underline{E}\Gamma,\kK))$ is an isomorphism for any finite subgroup $F$ of $\Gamma$: this is true since $\underline{E}\Gamma$ is $F$-uniformly contractible.
\pfend

\begin{rem}
If $Z$ is a classifying space for proper actions, then there exists a $\Gamma$-equivariant homotopy $f:[0,1]\times Z\times Z\to Z$ between the two projections $Z\times Z\to Z$. Then for every finite subgroup $F$ of $Z$ and every $F$-fixed point $a\in Z$, the map $(t,z)\mapsto f(t,z,a)$ is a $F$-equivariant homotopy between the identity and a constant map, but it is not necessarily uniformly continuous. This explains the extra condition that $\underline{E}\Gamma$ is $F$-uniformly contractible. On the other hand, it would be surprising if there existed a group for which $\underline{E}\Gamma$ is $\Gamma$-compact but not uniformly contractible.
\end{rem}

\section{Final remarks}
One of the main advantages of the coarse category is that it is much more flexible than the category of discrete groups. For instance, the coarse Baum--Connes map is invariant under coarse homotopy equivalence. It is natural to ask whether the (full or reduced) coarse Baum--Connes map with coefficients is also invariant under coarse homotopy equivalence, but the answer is probably not obvious.
To see why, let us consider a coarse map $f:X\to Y$. Let $Z=X\amalg Y$, endowed with the largest distance $d$ such that $d(x,f(x))=1$ $\forall x\in X$, $d_{|Y\times Y}\le d_Y$, $d_{|X\times X}\le d_X$. Then $Z$ is coarsely equivalent to $Y$, so there are coarse maps
$$C^*(X)\to C^*(Z)\leftarrow C^*(Y)$$
where the second map $C^*(Z)\leftarrow C^*(Y)$ is a Morita equivalence. Therefore, $f$ induces a map $K(C^*(X))\to K(C^*(Y))$.

To generalize such a construction to the conjecture with coefficients, it would be natural to expect similar maps on the groupoid level $G(X)\to G(Z)\leftarrow G(Y)$. However, the natural inclusion $X\to Z$ generally does not induce a map $G(X)\to G(Z)$ (unless $X\to Y$ is a coarse embedding). It does induce a map $G'(X)\to G'(Z)$, but the natural inclusion $G'(Y)\to G'(Z)$ is not a Morita equivalence.

\end{document}